# Order determination in general vector autoregressions

## Bent Nielsen[1]

*University of Oxford*

**Abstract:** In the application of autoregressive models the order of the model is often estimated using either a sequence of likelihood ratio tests, a likelihood based information criterion, or a residual based test. The properties of such procedures has been discussed extensively under the assumption that the characteristic roots of the autoregression are stationary. While non-stationary situations have also been considered the results in the literature depend on conditions to the characteristic roots. It is here shown that these methods for lag length determination can be used regardless of the assumption to the characteristic roots and also in the presence of deterministic terms. The proofs are based on methods developed by C. Z. Wei in his joint work with T. L. Lai.

## 1. Introduction

Order determination for stationary autoregressive time series has been discussed extensively in the literature. The three prevailing methods are either to test redundancy of the last lag using a likelihood based test, to estimate the lag length consistently using an information criteria, or to investigate the residuals of a fitted model with respect to autocorrelation. It is shown that these methods can be used regardless of any assumptions to the characteristic roots. This is important in applications, as the question of lag length can be addressed without having to locate the characteristic roots.

The statistical model is given by a $p$-dimensional time series $X_t$ of length $K + T$ satisfying a $K$th order vector autoregressive equation

$$(1.1) \qquad X_t = \sum_{l=1}^{K} A_l X_{t-l} + \mu D_t + \varepsilon_t, \qquad t = 1, \ldots, T,$$

conditional on the initial values $X_0, \ldots, X_{1-K}$. The effective sample will remain $X_1, \ldots, X_T$ when discussing autoregressions with $k < K$ to allow comparison of likelihood values. The component $D_t$ is a vector of deterministic terms such as a constant, a linear trend, or seasonal dummies. For the sake of defining a likelihood function it is initially assumed that the innovations, $(\varepsilon_t)$, are independently, identically normal, $\mathsf{N}_p(0, \Omega)$, distributed and independent of the initial values.

The aim is to determine the largest non-trivial order for the time series, $k_0$ say with $0 \le k_0 \le K$, so $A_{k_0} \ne 0$ and $A_j = 0$ for $j > k_0$. Three approaches are available of which the first is based on a likelihood ratio test for $A_k = 0$ where $1 \le k \le K$. The log likelihood ratio test statistic is

$$\mathsf{LR}\,(A_k = 0) = T \log \det \hat{\Omega}_{k-1} - T \log \det \hat{\Omega}_k,$$

[1]Department of Economics, University of Oxford & Nuffield College, Oxford OX1 1NF, UK, e-mail: bent.nielsen@nuffield.ox.ac.uk

*AMS 2000 subject classifications:* primary 62M10; secondary 62F10.

*Keywords and phrases:* autoregression, lag length, information criteria.

where $\hat{\Omega}_{k-j}$ is the conditional maximum likelihood estimator based on the observations $X_1, \ldots, X_T$ given the initial values, see (3.2) below. The statistic $\mathsf{LR}$ is proved to be asymptotically $\chi^2$ under the hypothesis $k_0 < k$, generalising results for the purely non-explosive case. Since the result does not depend on the characteristic roots, it can be used for lag length determination before locating the characteristic roots.

The second approach is to estimate $k_0$ by the argument $\hat{k}$ that maximises a penalised likelihood, or equivalently, minimises an information criteria of the type

$$(1.2) \qquad \Phi_j = \log \det \hat{\Omega}_j + j \frac{f(T)}{T}, \qquad j = 0, \ldots, K.$$

In the literature there are several candidates for the penalty function $f$. Akaike has $f(T) = 2p^2$, Schwarz [23] has $f(T) = p^2 \log T$ while Hannan and Quinn [10] and Quinn [22] have $f(T) = 2p^2 \log \log T$. For stationary processes without deterministic components it has been shown that the estimator $\hat{k}$ is weakly consistent if $f(T) = \mathrm{o}(T)$ and $f(T) \to \infty$ as $T$ increases, while Hannan and Quinn show, for $p = 1$, that strong consistency is obtained if $f(T) = \mathrm{o}(T)$ and $\liminf_{T\to\infty} f(T)/\log \log T > 2$, while strong consistency cannot be obtained if $\limsup_{T\to\infty} f(T)/\log \log T < 2$. In other words the estimators of Hannan and Quinn and of Schwarz are consistent while Akaike's estimator is inconsistent. Some generalisations to non-explosive processes have been given by for instance Paulsen [20], Pötscher [21] and Tsay [24]. Pötscher also considered the purely explosive case but did not obtain a common feasible rate for $f(T)$ for the explosive and the non-explosive case. In the following consistency is shown for a penalty function $f(T)$ not depending on the characteristic roots, showing that the penalised likelihood approach also can be applied to lag length determination prior to locating the characteristic roots.

A third approach is a residual based mis-specification test. This is implemented in particular in econometric computer packages. In a first step the residuals, $\hat{\varepsilon}_t$ say, are computed from the model (1.1) with $k - 1$ lags, say. In a second step an auxillary regression is considered where $\hat{\varepsilon}_t$ is regressed on lagged values as well as the regressors in equation (1.1). It is argued that a test based on the squared multiple correlation arising from the auxillary regression is asymptotically equivalent to the above mentioned likelihood ratio test statistic also in the general case.

Like the work of Pötscher [21] the proofs in this paper are based on the joint work of C. Z. Wei and T. L. Lai on the strong consistency of least squares estimators in autoregressions, see for instance Lai and Wei [15]. As pointed out in Pötscher's Remark 1 to his Theorem 3.3 these results are not quite strong enough to facilitate common feasible rates for the penalty function. Two important ingredients in the presented proofs are therefore an algebraic decomposition exploiting partitioned inversion along with a generalisation of Lai and Wei's work given by Nielsen [17]. Whereas the former paper is concerned with showing that the least squares estimator for the autoregressive estimator is consistent, the latter paper provides a more detailed discussion of the rate of consistency as well as it allows deterministic terms in the autoregression.

The following notation is used throughout the paper: For a quadratic matrix $\alpha$ let $\mathrm{tr}(\alpha)$ denote the trace and $\lambda(\alpha)$ the set of eigenvalues, so that $|\lambda(\alpha)| < 1$ means that all eigenvalues have absolute value less than one. When $\alpha$ is also symmetric then $\lambda_{\min}(\alpha)$ and $\lambda_{\max}(\alpha)$ denote the smallest and the largest eigenvalue respectively. The abbreviations $a.s.$ and $\mathsf{P}$ are used for properties holding almost surely and in probability, respectively.

## 2. Results

Before presenting the results the assumptions and notation is set up. Then the results follow for the three approaches.

### 2.1. Assumptions and notation

The asymptotic analysis is to a large extent based on results of Lai and Wei [15] with appropriate modifications to the situation with deterministic terms in Nielsen [17]. Following that analysis the assumption to the innovations of independence and normality made above can be relaxed so that the sequence of innovations $(\varepsilon_t)$ is a martingale difference sequence with respect to an increasing sequence of $\sigma$-fields $(\mathcal{F}_t)$, that is: the innovations $X_{1-k}, \ldots, X_0$ are $\mathcal{F}_0$-measurable and $\varepsilon_t$ is $\mathcal{F}_t$-measurable with $\mathsf{E}(\varepsilon_t | \mathcal{F}_{t-1}) \overset{a.s.}{=} 0$, which is assumed to satisfy

$$(2.1) \qquad \sup_t \mathsf{E}\{(\varepsilon_t'\varepsilon_t)^{\lambda/2} | \mathcal{F}_{t-1}\} \overset{a.s.}{<} \infty \qquad \text{for some } \lambda > 4.$$

To establish an asymptotic theory for the LR-statistic it is assumed that

$$(2.2) \qquad \mathsf{E}(\varepsilon_t \varepsilon_t' | \mathcal{F}_{t-1}) \overset{a.s.}{=} \Omega,$$

where $\Omega$ is positive definite. For the asymptotic theory for the information criteria this can be relaxed to

$$(2.3) \qquad \liminf_{t \to \infty} \lambda_{\min} \mathsf{E}(\varepsilon_t \varepsilon_t' | \mathcal{F}_{t-1}) \overset{a.s.}{>} 0.$$

The deterministic term $D_t$ is a vector of terms such as a constant, a linear trend, or periodic functions like seasonal dummies. Inspired by Johansen [13] the deterministic terms are required to satisfy the difference equation

$$(2.4) \qquad D_t = \mathbf{D} D_{t-1},$$

where $\mathbf{D}$ has characteristic roots on the complex unit circle. For example,

$$\mathbf{D} = \begin{pmatrix} 1 & 0 \\ 1 & -1 \end{pmatrix} \qquad \text{with} \qquad D_0 = \begin{pmatrix} 1 \\ 1 \end{pmatrix}$$

will generate a constant and a dummy for a biannual frequency. The deterministic term $D_t$ is assumed to have linearly independent coordinates. That is:

$$(2.5) \qquad |\lambda(\mathbf{D})| = 1, \qquad \operatorname{rank}(D_1, \ldots, D_{\dim \mathbf{D}}) = \dim \mathbf{D}.$$

In the analysis it is convenient to introduce the companion form

$$(2.6) \qquad \begin{pmatrix} \mathbf{X}_t \\ D_t \end{pmatrix} = \begin{pmatrix} \mathbf{B} & \boldsymbol{\mu} \\ 0 & \mathbf{D} \end{pmatrix} \begin{pmatrix} \mathbf{X}_{t-1} \\ D_{t-1} \end{pmatrix} + \begin{pmatrix} \mathbf{e}_t \\ 0 \end{pmatrix},$$

where $\mathbf{X}_{t-1} = (X_{t-1}', \ldots, X_{t-k+1}')'$ and

$$\mathbf{B} = \left\{ \begin{array}{ccc} A_1 \cdots A_{k-2} & A_{k-1} \\ I_{p(k-2)} & 0 \end{array} \right\}, \quad \boldsymbol{\iota} = \left\{ \begin{array}{c} I_p \\ 0_{(k-2)p \times p} \end{array} \right\}, \quad \boldsymbol{\mu} = \boldsymbol{\iota} \mu \mathbf{D}, \quad \mathbf{e}_t = \boldsymbol{\iota} \varepsilon_t.$$

The process $\mathbf{X}_t$ can be decomposed using a similarity transformation. Following Herstein ([11], p. 308) there exists a regular, real matrix $M$ that block-diagonalises $\mathbf{B}$ so that $M\mathbf{B}M^{-1} = \mathrm{diag}(\mathbf{U}, \mathbf{V}, \mathbf{W})$ is a real block diagonal matrix where the eigenvalues of the diagonal blocks $\mathbf{U}, \mathbf{V}, \mathbf{W}$ satisfy $|\lambda(\mathbf{U})| < 1$, $|\lambda(\mathbf{V})| = 1$, and $|\lambda(\mathbf{W})| > 1$. Any of the blocks $\mathbf{U}, \mathbf{V}, \mathbf{W}$ can be empty matrices, so if for instance $|\lambda(\mathbf{B})| < 1$ then $\mathbf{U} = \mathbf{B}$ and $\dim \mathbf{V} = \dim \mathbf{W} = 0$. The process $\mathbf{X}_t$ can therefore be decomposed as

$$(2.7) \qquad M\mathbf{X}_t = \begin{pmatrix} U_t \\ V_t \\ W_t \end{pmatrix} = \begin{pmatrix} \mathbf{U} & 0 & 0 & \mu_U \\ 0 & \mathbf{V} & 0 & \mu_V \\ 0 & 0 & \mathbf{W} & \mu_W \end{pmatrix} \begin{pmatrix} U_{t-1} \\ V_{t-1} \\ W_{t-1} \\ D_t \end{pmatrix} + \begin{pmatrix} e_{U,t} \\ e_{V,t} \\ e_{W,t} \end{pmatrix}.$$

Finally, there exists a constant $\tilde{\mu}_U$, see Nielsen ([17], Lemma 2.1), so

$$(2.8) \qquad U_t = \tilde{U}_t + \tilde{\mu}_U D_t \qquad \text{where} \qquad \tilde{U}_t = \mathbf{U}\tilde{U}_{t-1} + e_{U,t}.$$

### 2.2. Likelihood ratio test statistics

The likelihood ratio test statistic is known to be asymptotically $\chi^2$ in the stationary case where $|\lambda(\mathbf{B})| < 1$ and $\mathbf{D} = 1$, see Lütkepohl ([16], Section 4.2.2). Here the result is shown to hold regardless of the assumptions to $\mathbf{B}$ and $\mathbf{D}$. Thus, the likelihood ratio test can be used before locating the charateristic roots.

**Theorem 2.1.** *Suppose Assumptions* (2.1), (2.2), (2.5) *are satisfied and* $k_0 < k$. *Then* $\mathsf{LR}(A_k = 0)$ *is asymptotically* $\chi^2(p^2)$.

Since the likelihood ratio test statistic is based on partial correlations it follows from Theorem 2.1 that partial correlograms that are computed from partial correlograms can be used regardless of the location of the characteristic roots. Often correlograms are, however, based on the Yule-Walker estimators, which assume stationarity. For non-stationary autoregressions that can lead to misleading inference. Nielsen [18] provides a more detailed discussion.

**Remark 2.2.** The fourth order moment condition, $\lambda > 4$, in Assumption (2.1) is used twice in the proof. First, to ensure that the residuals from regressing $\varepsilon_t$ on the explosive term $W_{t-1}$ do not depend asymptotically on $W_{t-1}$. As discussed in Remark 3.7 it suffices that $\lambda > 2$ if either of the following conditions hold:

(I,a)  $\dim \mathbf{W} = 0$.

(I,b)  $\dim \mathbf{W} > 0$ and $\varepsilon_t$ independent, identically distributed.

Secondly, to ensure that $\varepsilon_t \varepsilon_{t-1}$ has second moments when applying a Central Limit Theorem. As discussed in Remark 3.12, it suffice that $\lambda > 2$ if

(II)  the innovations $\varepsilon_t$ are independent.

The test statistic considered above is for a hypothesis concerning a single lag. This can be generalised to a hypothesis concerning several lags, $m$ say, where $k + m - 1 \leq K$.

**Theorem 2.3.** *Suppose Assumptions* (2.1), (2.2), (2.5) *are satisfied and* $k_0 < k$. *Then* $\mathsf{LR}(A_k = \cdots = A_{k+m-1} = 0)$ *is asymptotically* $\chi^2(p^2 m)$.

### 2.3. Information criteria

The next two results concern consistency of a lag length estimator arising from use of information criterions. The proof has two distinct parts. First, it is argued that the lag length estimator $\hat{k}$ is not under-estimating, and, secondly, that it is not over-estimating. The first part is the easy one to establish. This result holds for all of the penalty functions discussed in the introduction under weak conditions to the innovations.

**Theorem 2.4.** *Suppose Assumptions (2.1), (2.3), (2.5) are satisfied with $\lambda > 2$ only and $f(T) = \mathrm{o}(T)$. Then $\liminf_{T\to\infty} \hat{k} \overset{a.s.}{\geq} k_0$.*

This result has previously been established in the univariate case without deterministic terms so $p = \dim X = 1$ and $\dim \mathbf{D} = 0$ by Pötscher (1989, Theorem 3.3). For the purely explosive case $|\lambda(\mathbf{B})| > 1$ his Theorem 3.2 shows the above result under the weaker condition $f(T) = \mathrm{o}(T^2)$. A version holding in probabilty has been shown for the non-explosive case $|\lambda(\mathbf{B})| \leq 1$ and $\mathbf{D} = 1$ by Paulsen [20] and Tsay [24].

Results showing that the lag length is not overestimating are harder to establish. Various weak and strong results can be obtained depending on the number of conditions that are imposed.

**Theorem 2.5.** *Suppose Assumptions (2.1), (2.5) are satisfied. Then*

(i) *If $f(T) \to \infty$ and Assumption (2.2) holds then $\mathsf{P}(\hat{k} \leq k_0) \to 1$.*

(ii) *If $f(T)/\log T \to \infty$ and Assumption (2.3) holds then $\limsup_{T\to\infty} \hat{k} \overset{a.s.}{\leq} k_0$.*

(iii) *If $f(T)/\{(\log\log T)^{1/2}(\log T)^{1/2}\} \to \infty$, Assumption (2.3) holds, and the parameters satisfy the condition $(A)$ that $\mathbf{V}$ and $\mathbf{D}$ have no common eigenvalues then $\limsup_{T\to\infty} \hat{k} \overset{a.s.}{\leq} k_0$.*

(iv) *If $f(T)/\log\log T \to \infty$, Assumption (2.3) holds, and either $(B) \dim \mathbf{D} = 0$ with $\mathbf{V} = 1$ or $(C) \dim \mathbf{V} = 0$ then $\limsup_{T\to\infty} \hat{k} \overset{a.s.}{\leq} k_0$.*

(v) *Suppose Assumption (2.2) holds, and either $(B)$ or $(C)$ holds then*

    (a) *If $\liminf_{T\to\infty} (2\log\log T)^{-1} f(T) \overset{a.s.}{>} p^2$ then $\limsup_{T\to\infty} \hat{k} \overset{a.s.}{\leq} k_0$.*

    (b) *If $\limsup_{T\to\infty} (2\log\log T)^{-1} f(T) \overset{a.s.}{<} 1$ then $\hat{k} \overset{a.s.}{\nrightarrow} k_0$.*

By combining Theorems 2.4, 2.5 consistency results can be obtained. For instance Theorem 2.4 in combination with Theorem 2.5(i) shows that $\hat{k} \overset{\mathsf{P}}{\to} k_0$ if the penalty function satisfies $f(T) \to \infty$ and $f(T) = \mathrm{o}(T)$. This includes Hannan and Quinn's and Schwarz's penalty functions, but excludes that of Akaike as usually found. Likewise, Theorem 2.4 in combination with Theorem 2.5(ii) show that $\hat{k} \overset{a.s.}{\to} k_0$ if the penalty function satisfies $f(T)/\log T \to \infty$ and $f(T) = \mathrm{o}(T)$. These results are the first to present conditions to the penalty function ensuring consistency that are not depending on the parameter $\mathbf{B}$ and $\mathbf{D}$. This implies that the information criteria can be used before locating the charateristic roots.

It remains an open problem, however, to establish strong consistency of the Schwarz and the Hannan-Quinn estimators for general values of $\mathbf{V}$ and $\mathbf{D}$. Theorem 2.4 in combination with Theorem 2.5(iii) shows that the Schwarz estimator is strongly consistent when $(A)$ holds so $\mathbf{V}$ and $\mathbf{D}$ have no common eigenvalues. Theorem 2.4 combined with Theorem 2.5(v) shows that the Hannan-Quinn estimator is strongly consistent when either $(B) \dim \mathbf{D} = 0$ with $\mathbf{V} = 1$ or $(C) \dim \mathbf{V} = 0$ holds. This is the first strong consistency result for the Hannan-Quinn estimator in the nonstationary case.

**Remark 2.6.** In Theorem 2.5 the fourth order moment condition $\lambda > 4$ in Assumption (2.1) can be relaxed to $\lambda > 2$ under certain condions to the parameters. Recall the conditions stated in Remark 2.2 which are

(I,a)  dim $\mathbf{W} = 0$.
(I,b)  dim $\mathbf{W} > 0$ and $\varepsilon_t$ independent, identically distributed.
 (II)  the innovations $\varepsilon_t$ are independent.

As discussed in Remark 3.13 it holds:
Result (i) can be relaxed if (II) holds along with either (I,a) or (I,b).
Results (ii), (iii), (iv) can be relaxed if (I,a) holds.
Result (v) cannot be relaxed with the present proof.

A number of related results are available in the literature.

The weak consistency results in (i) has been shown for the non-explosive case $|\lambda(\mathbf{B})| \leq 1$ and $\mathbf{D} = 1$ by Paulsen [20] and Tsay [24].

The $(\log \log T)^{1/2}(\log T)^{1/2}$ rate discussed in Theorem 2.5(iii) and Remark 2.6(iii) is an improvement over the $\log T$ rates discussed by for instance Pötscher [21] and Wei [25]. These authors discuss the univariate case without deterministic terms so $p = \dim X = 1$ and $\dim \mathbf{D} = 0$, in which case $\mathbf{V}$ and $\mathbf{D}$ trivially have no common eigenvalues. *First*, Pötscher ([21], Theorem 3.1) shows an under-estimation result for rates satisfying $f(T)/\log T \to \infty$ in the non-explosive case so $|\lambda(\mathbf{B})| \leq 1$, hence $\dim \mathbf{W} = 0$, but with Assumption (2.3) replaced by the weaker condition that $\liminf_{T \to \infty} T^{-1} \sum_{t=1}^{T} \mathsf{E}(\varepsilon_t^2 | \mathcal{F}_{t-1}) \overset{a.s.}{\to} 0$. Pötcher's Theorem 3.2 concerning under-estimation in the purely explosive case so $|\lambda(\mathbf{B})| > 1$ requires $\liminf_{T \to \infty} f(T)/T > 0$ *a.s.* with just $\lambda > 2$ in Assumption (2.1). The Remark 1 to his Theorem 3.3 points out that his results do not provide a common feasibility rate for autoregressions with both explosive and non-explosive roots in that $f(T) = \mathrm{o}(T)$ is required for the over-estimation result, whereas $\liminf_{T \to \infty} f(T)/T > 0$ *a.s.* is required for the under-estimation result. *Secondly*, Theorem 3.6 of Wei [25] goes a step further in showing the over-estimation result for the rate $f(T) = \log T$ for the non-explosive case so $\dim \mathbf{W} = 0$.

The optimal $\log \log T$ rates in (v) were originally suggested by Hannan and Quinn [10] and Quinn [22] for the case where $|\lambda(\mathbf{B})| < 1$, $\dim \mathbf{D} = 0$. A full generalisation cannot be made at present as the proof hinges on proving that the smallest eigenvalue of the average of the squared residual from regressing $V_{t-1}$ on $D_t$, that is $T^{-1-\eta} \sum_{t=1}^{T} (V_{t-1}|D_t)(V_{t-1}|D_t)'$, has positive limit points for some $\eta > 0$. This result can only be established in two special cases: first, if $\dim \mathbf{V} = 0$ the issue is irrelevant, and secondly, if $\mathbf{V} = 1$ and $\dim \mathbf{D} = 0$ this follows from the law of iterated logarithms by Donsker and Varadhan [6]. A more detailed discussion is given in Lemma 3.5(iv) in the Appendix.

The strong $\log \log T$ rate in Theorem 2.5(iv) and Remark 2.6(iv) has previously been established in the purely stable, univariate case without deterministic terms, so $p = \dim X = 1$ and $\dim \mathbf{D} = 0$ and $|\lambda(\mathbf{B})| < 1$, and hence $\dim \mathbf{W} = 0$, see Pötscher ([21], Theorem 3.4). Once again, his result only requires $\liminf_{T \to \infty} T^{-1} \sum_{t=1}^{T} \mathsf{E}(\varepsilon_t^2 | \mathcal{F}_{t-1}) \to 0$ *a.s.* instead of Assumption 2.3.

## *2.4. Residual based mis-specification testing*

The third approach is to fit the model (1.1) with $k-1$ lags and analyse the residuals for autocorrelation of order up to $m$. The maximal lag length parameter $K$ is here required to be at least $k - 1$. This is done in two steps. First the residuals $\hat{\varepsilon}_t$ are found for the regression (1.1) with $t = 1, \ldots, T$ and $k - 1$ lags. In the second step $\hat{\varepsilon}_t$ is analysed in an auxiliary regression for $t = m + 1, \ldots, T$, where $\hat{\varepsilon}_t$ is regressed on $\hat{\varepsilon}_{t-1}, \ldots, \hat{\varepsilon}_{t-m}$ as well as the original regressors $\mathbf{X}_{t-1} = (X'_{t-1}, \ldots, X'_{t-k+1})'$ and $D_t$. The original regressors are included to mimic the above likelihood analysis where $\mathbf{X}_{t-1}, D_t$ are partialled out from $X_t$ and $X_{t-k}$. A test based on the squared sample correlation of the variables in the auxiliary regression is asymptotically equivalent to the likelihood ratio tests, so the degrees of freedom do not include the dimension of $\mathbf{X}_{t-1}, D_t$. In the multivariate case, $p > 1$, the test can be implemented in three ways, using either a simultaneous test, a marginal test or a conditional test.

The joint test, is based on the test statistic $\text{tr}(TR^2)$, where $R^2$ is the squared sample multiple correlation of $\hat{\varepsilon}_t$ and $(\hat{\varepsilon}'_{t-1}, \ldots, \hat{\varepsilon}'_{t-m}, \mathbf{X}'_{t-1}, D'_t)'$.

The other two tests are based on a $q$-dimensional subset of the $p$ components of $\varepsilon_t$. As the equations in the model equation (1.1) can be permuted there is no loss of generality in focussing on the first $q$ components. Thus, partition

$$\varepsilon_t = \begin{pmatrix} \varepsilon_{t,1} \\ \varepsilon_{t,2} \end{pmatrix}, \qquad X_t = \begin{pmatrix} X_{t,1} \\ X_{t,2} \end{pmatrix},$$

where $\varepsilon_{t,1}$ and $X_{t,1}$ are $q$-dimensional.

The marginal model consists of the first $q$ equations of (1.1), that is $X_{t,1}$ given $\mathbf{X}_{t-1}, D_t$. The marginal test is then based on the squared sample multiple correlation, $R^2_{\text{marg}}$ say, of $\hat{\varepsilon}_{t,1}$ and $(\hat{\varepsilon}'_{t-1,1}, \ldots, \hat{\varepsilon}'_{t-m,1}, \mathbf{X}'_{t-1}, D'_t)$.

The conditional model consists of the first $q$ equations of (1.1) given $X_{t,2}$, that is $X_{t,1}$ given $X_{t,2}, \mathbf{X}_{t-1}, D_t$. The conditional test is based on the squared sample multiple correlation, $R^2_{\text{cond}}$ say, of $\hat{\varepsilon}_{t,1}$ and $(\hat{\varepsilon}'_{t-1,1}, \ldots, \hat{\varepsilon}'_{t-m,1}, X'_{t,2}, \mathbf{X}'_{t-1}, D'_t)$.

The following asymptotic result can be established.

**Theorem 2.7.** *Suppose Assumptions (2.1), (2.2), (2.5) are satisfied and $k_0 < k$. Then $\text{tr}(TR^2)$ is asymptotically $\chi^2(p^2 m)$, while $\text{tr}(TR^2_{marg})$ and $\text{tr}(TR^2_{cond})$ are asymptotically $\chi^2(q^2 m)$.*

Sometimes these test are implemented so that the auxiliary regression is carried out for $t = 1, \ldots, T$ rather than $t = m+1, \ldots, T$ with the convention that $\hat{\varepsilon}_0 = \cdots = \hat{\varepsilon}_{1-m} = 0$. Variants of the tests have been considered, in particular for the univariate case, by Durbin [7], Godfrey [8], Breusch [3] and Pagan [19]. Those variants have been argued to be score/Lagrange multiplier type tests and asymptotic theory has been established for the stationary case $|\lambda(\mathbf{B})| < 1$.

## 3. Proofs

The likelihood ratio test statistic for testing $A_k = 0$ is given by

$$
\begin{aligned}
LR\left(A_k = 0\right) &= -T \log \det(\hat{\Omega}_{k-1}^{-1} \hat{\Omega}_k) \\
&= -T \log \det\{I_p - \hat{\Omega}_{k-1}^{-1}(\hat{\Omega}_{k-1} - \hat{\Omega}_k)\},
\end{aligned}
\tag{3.1}
$$

where $\hat{\Omega}_k$ and $\hat{\Omega}_{k-1}$ represent the unrestricted and restricted maximum likelihood estimators for the variance matrix defined below. In the following first some notation is introduced. Then comes an asymptotic analysis of $\hat{\Omega}_{k-1}$ and $\hat{\Omega}_{k-1} - \hat{\Omega}_k$ and finally proofs of the main theorems follow.

### 3.1. Notation

It is convenient to introduce some notation to handle $\hat{\Omega}_{k-1}$ as well as $\hat{\Omega}_{k-1} - \hat{\Omega}_k$. Thus, let the residuals from the partial regressions of $X_t$ and $X_{t-k}$ on $\mathbf{X}_{t-1} = (X'_{t-1}, \ldots, X'_{t-k+1})'$ and the deterministic components $D_t$ be denoted

$$(X_t | \mathbf{X}_{t-1}, D_t), \qquad (X_{t-k} | \mathbf{X}_{t-1}, D_t).$$

When the hypothesis, $A_k = 0$, is satisfied then $(X_t | \mathbf{X}_{t-1}, D_t) = (\varepsilon_t | \mathbf{X}_{t-1}, D_t)$ and therefore the restricted variance estimator is given by

$$(3.2) \qquad \hat{\Omega}_{k-1} = \frac{1}{T} \sum_{t=1}^{T} (\varepsilon_t | \mathbf{X}_{t-1}, D_t) \, (\varepsilon_t | \mathbf{X}_{t-1}, D_t)'.$$

Most of the analysis in the proof relates to $\hat{\Omega}_{k-1} - \hat{\Omega}_k$ so it is helpful to define

$$Q(Z_t) = \sum_{t=1}^{T} \varepsilon_t Z'_t \left( \sum_{t=1}^{T} Z_t Z'_t \right)^{-1} \sum_{t=1}^{T} Z_t \varepsilon'_t,$$

for any time series $Z_t$. It follows that $T(\hat{\Omega}_{k-1} - \hat{\Omega}_k) = Q(X_{t-k} | \mathbf{X}_{t-1}, D_t)$. Occasionally the following notation will be used: For a matrix $\alpha$ let $\alpha^{\otimes 2} = \alpha \alpha'$.

### 3.2. Asymptotic analysis of $\hat{\Omega}_{k-1}$

Asymptotic expressions for the restricted least squares variance estimator $\hat{\Omega}_{k-1}$ are given by Nielsen ([17], Corollary 2.6, Theorem 2.8):

**Lemma 3.1.** *Suppose $A_k = 0$ and that the Assumptions $(2.1), (2.3), (2.5)$ are satisfied with $\lambda > 2$. Then, for all $\xi < 1 - 2/\lambda$ it holds*

$$\hat{\Omega}_{k-1} \overset{a.s.}{=} \frac{1}{T} \sum_{t=1}^{T} \varepsilon_t \varepsilon'_t + \mathrm{o}(T^{-\xi}),$$

*If in addition Assumption $(2.2)$ is satisfied then for all $\zeta < \min(\xi, 1/2)$ it holds*

$$\hat{\Omega}_{k-1} \overset{a.s.}{=} \Omega + \mathrm{o}(T^{-\zeta}).$$

### 3.3. Asymptotic analysis of $\hat{\Omega}_{k-1} - \hat{\Omega}_k$

The analysis of the term $\hat{\Omega}_{k-1} - \hat{\Omega}_k$ is specific to the order selection problem. For the sake of finding the asymptotic distribution of the likelihood ratio test statistic the aim is to express $\hat{\Omega}_{k-1} - \hat{\Omega}_k$ in terms of a stationary process $Y_t$ as

$$(3.3) \qquad T(\hat{\Omega}_{k-1} - \hat{\Omega}_k) = Q(X_{t-k} | \mathbf{X}_{t-1}, D_t) = Q(Y_{t-1}) + \mathrm{o}_{\mathsf{P}}(1),$$

which in turn can be proved to be asymptotically $\chi^2$ by a Central Limit Theorem. The result (3.3) reduces trivially to an equality with $Y_{t-1} = \varepsilon_{t-1}$ when testing $A_1 = 0$, so only the case $k > 1$ will need consideration in the remainder of this subsection. On the way to prove the above result some related expressions holding under weaker assumptions emerge which can be used for proving the consistency results for the estimator of the lag length, $\hat{k}$.

In the following $\hat{\Omega}_{k-1} - \hat{\Omega}_k$ is first decomposed into seven terms. It is then shown that the three leading term can be written as $Q(Y_{t-1})$ as in (3.3) and that the remaining four terms are asymptotically vanishing.

### 3.3.1. *Decomposition of* $\hat{\Omega}_{k-1} - \hat{\Omega}_k$

The first decomposition is a purely algebraic result based on the formula for partitioned inversion.

**Lemma 3.2.** *Suppose* $A_k = 0$. *Then it holds*

$$Q\left(X_{t-k}|\mathbf{X}_{t-1}, D_t\right) = Q\left(\mathbf{X}_{t-2}|D_t\right) - Q\left(\mathbf{X}_{t-1}|D_t\right) + Q\left(\varepsilon_{t-1}|\mathbf{X}_{t-2}, D_t\right).$$

*Proof of Lemma 3.2.* By the formula for partitioned inversion it holds

$$(3.4) \qquad Q\left(\begin{array}{c}\mathbf{X}_{t-1}\\ X_{t-k}\end{array}\bigg| D_t\right) = Q\left(X_{t-k}|\mathbf{X}_{t-1}, D_t\right) + Q\left(\mathbf{X}_{t-1}|D_t\right),$$

of which $T(\hat{\Omega}_{k-1} - \hat{\Omega}_k) = Q\left(X_{t-k}|\mathbf{X}_{t-1}, D_t\right)$ is the first term on the left. Noting that $(\mathbf{X}'_{t-1}, X'_{t-k})' = (X'_{t-1}, \mathbf{X}'_{t-2})'$ a repeated use of the formula for partitioned inversion shows

$$(3.5) \quad Q\left(\begin{array}{c}\mathbf{X}_{t-1}\\ X_{t-k}\end{array}\bigg| D_t\right) = Q\left(\begin{array}{c}X_{t-1}\\ \mathbf{X}_{t-2}\end{array}\bigg| D_t\right) = Q\left(X_{t-1}|\mathbf{X}_{t-2}, D_t\right) + Q\left(\mathbf{X}_{t-2}|D_t\right).$$

Due to the model equation (1.1) with $A_k = 0$ and the property $D_t = \mathbf{D}D_{t-1}$ it follows $(X_{t-1}|\mathbf{X}_{t-2}, D_t) = (\varepsilon_{t-1}|\mathbf{X}_{t-2}, D_t)$. The desired expression then arise by rearranging the above expressions. $\square$

Asymptotic arguments are now needed. These arguments rely on Nielsen [17] which in turn represents a generalisation of the arguments of Lai and Wei [15]. The second step is therefore an asymptotic decomposition of the first two terms in Lemma 3.2 using that the processes $U_t, V_t, W_t$ are asymptotically uncorrelated.

**Lemma 3.3.** *Suppose* $A_k = 0$ *and that the Assumptions* (2.1), (2.3), (2.5) *are satisfied with* $\lambda > 2$. *Then, for* $j = 1, 2$,

$$(3.6) \qquad Q\left(\mathbf{X}_{t-j}|D_t\right) \stackrel{a.s.}{=} Q\left(U_{t-j}|D_t\right) + Q\left(V_{t-j}|D_t\right) + Q\left(W_{t-j}|D_t\right) + \mathrm{o}\left(1\right).$$

*Proof of Lemma 3.3.* Since $M\mathbf{X}_t = (U_t, V_t, W_t)$, see (2.7), it suffices to argue that the processes $U_t, V_t$ and $W_t$ are asymptotically uncorrelated so that the off-diagonal elements of $\sum_{t=1}^{T}(\mathbf{X}_{t-j}|D_t)(\mathbf{X}_{t-j}|D_t)'$ can be ignored in the asymptotic argument. This follows from Nielsen ([17], Theorem 6.4, 9.1, 9.2, 9.4), see also the summary in Table 2 of that paper. $\square$

### 3.3.2. *Eliminating explosive terms and regressors in stationary terms*

In combination Lemmas 3.2, 3.3 show that

$$T(\hat{\Omega}_{k-1} - \hat{\Omega}_k) \stackrel{a.s.}{=} Q\left(\varepsilon_{t-1}|\mathbf{X}_{t-2}, D_t\right) + Q\left(U_{t-2}|D_t\right) - Q\left(U_{t-1}|D_t\right)$$
$$+ Q\left(V_{t-2}|D_t\right) - Q\left(V_{t-1}|D_t\right) + Q\left(W_{t-2}|D_t\right) - Q\left(W_{t-1}|D_t\right) + \mathrm{o}\left(1\right).$$

Under mild conditions this can be reduced further so as to eliminate the terms involving the explosive component $W_t$ as well as the regressors in the terms involving the stationary component $U_t$.

**Lemma 3.4.** *Suppose $A_k = 0$ and that the Assumptions $(2.1)$, $(2.3)$, $(2.5)$ are satisfied, with $\lambda > 2$. Then,*

$$(3.7) \qquad T(\hat{\Omega}_{k-1} - \hat{\Omega}_k) \stackrel{a.s.}{=} Q\left(\varepsilon_{t-1}\right) + Q(\tilde{U}_{t-2}) - Q(\tilde{U}_{t-1}) + R_\varepsilon + R_V + \mathrm{o}\left(1\right),$$

*where*

$$(3.8) \quad R_\varepsilon = Q\left(\varepsilon_{t-1}|\mathbf{X}_{t-2}, D_t\right) - Q\left(\varepsilon_{t-1}\right), \quad R_V = Q\left(V_{t-2}|D_t\right) - Q\left(V_{t-1}|D_t\right).$$

*Proof of Lemma 3.4.* It suffices to prove, for $j = 1, 2$,

$$(3.9) \qquad\qquad\qquad Q\left(U_{t-j}|D_t\right) \stackrel{a.s.}{=} Q(\tilde{U}_{t-j}) + \mathrm{o}\left(1\right),$$

$$(3.10) \qquad Q\left(W_{t-2}|D_t\right) - Q\left(W_{t-1}|D_t\right) \stackrel{a.s.}{=} \mathrm{o}\left(1\right).$$

First, consider $(3.9)$. Because of $(2.8)$ then $(U_{t-j}|D_t) = (\tilde{U}_{t-j}|D_t)$. According to Nielsen ([17], Theorem 6.4) it holds for any $\eta > 0$ that

$$\left(\sum_{t=1}^{T} D_t D_t'\right)^{-1/2} \sum_{t=1}^{T} D_t \tilde{U}_{t-j}' \left(\sum_{t=1}^{T} \tilde{U}_{t-j}\tilde{U}_{t-j}'\right)^{-1/2} \stackrel{a.s.}{=} \mathrm{o}(T^{\eta-1/2}),$$

while Theorem 6.2 of the above paper shows $T^{-1}\sum_{t=1}^{T} \tilde{U}_{t-j}\tilde{U}_{t-j}'$ has positive definite limit points. This implies

$$\sum_{t=1}^{T} (\tilde{U}_{t-j}|D_t)(\tilde{U}_{t-j}|D_t)' \stackrel{a.s.}{=} \sum_{t=1}^{T} \tilde{U}_{t-j}\tilde{U}_{t-j}' \left\{1 + \mathrm{o}\left(T^{2\eta-1}\right)\right\}.$$

Theorem 2.4 of the above paper shows $\sum_{t=1}^{T} \varepsilon_t D_t'(\sum_{t=1}^{T} D_t D_t')^{-1/2} = \mathrm{o}(T^\eta)$ implying

$$\sum_{t=1}^{T} \varepsilon_t (\tilde{U}_{t-j}|D_t)' \stackrel{a.s.}{=} \sum_{t=1}^{T} \varepsilon_t \tilde{U}_{t-j}' + \mathrm{o}(T^{2\eta}).$$

That theorem also shows $\sum_{t=1}^{T} \varepsilon_t \tilde{U}_{t-j}'(\sum_{t=1}^{T} \tilde{U}_{t-j}\tilde{U}_{t-j}')^{-1/2} = \mathrm{o}(T^\eta)$. In combination these results show the desired result.

Secondly, consider $(3.10)$. Note first that $W_{t-1} = \mathbf{W}W_{t-2} + \mu_W D_{t-1} + e_{W,t-1}$ by $(2.7)$ while $D_{t-1} = \mathbf{D}^{-1}D_t$, implying $(W_{t-1}|D_t) = (\mathbf{W}W_{t-2} + e_{W,t-1}|D_t)$. This gives rise to the expansions

$$\sum_{t=1}^{T} (W_{t-1}|D_t)^{\otimes 2} = \sum_{t=1}^{T} (\mathbf{W}W_{t-2}|D_t)^{\otimes 2} \left(1 + f_T\right),$$

$$\sum_{t=1}^{T} (W_{t-1}|D_t)\varepsilon_t = \sum_{t=1}^{T} (\mathbf{W}W_{t-2}|D_t)\varepsilon_t + c_T,$$

where $f_T = \mathrm{O}(d_T^{-1/2}a_T) + d_T^{-1}b_T$ and

$$a_T = d_T^{-1/2}\sum_{t=1}^{T} (\mathbf{W}W_{t-2}|D_t)\,e_{W,t-1}, \qquad b_T = \sum_{t=1}^{T} (e_{W,t-1}|D_t)^{\otimes 2},$$

$$c_T = \sum_{t=1}^{T} (e_{W,t-1}|D_t)\varepsilon_t, \qquad\qquad d_T = \sum_{t=1}^{T} (\mathbf{W}W_{t-2}|D_t)^{\otimes 2}.$$

Using Nielsen ([17], Theorems 2.4, 6.2, 6.4) it is seen that

$$b_T \stackrel{a.s.}{=} \mathrm{O}(T), \quad c_T \stackrel{a.s.}{=} \mathrm{o}(T^{1/2+\eta}).$$

It follows from Nielsen ([17], Theorems 2.4, 9.1 and Corollary 7.2) that

$$Q\left(W_{t-j}|D_t\right) \stackrel{a.s.}{=} \mathrm{o}\left(T\right), \quad a_T \stackrel{a.s.}{=} \mathrm{o}(T^{1/2}), \quad d_T^{-1} \stackrel{a.s.}{=} \mathrm{o}\left(\rho^{-T}\right),$$

for some $\rho > 0$. This implies that $f_T$ is exponentially decreasing. The desired result follows by expanding $Q(W_{t-1}|D_t)$ in terms of $Q(W_{t-2}|D_t)$ as

$$\left[Q\left(W_{t-2}|D_t\right) + d_T^{-1/2}c_T\mathrm{O}\left\{Q\left(W_{t-2}|D_t\right)\right\}^{1/2} + c_T'd_T^{-1}c_T\right](1+f_T),$$

and using the established orders of magnitude. □

### 3.3.3. Eliminating unit root terms and regressors in innovation terms

The terms $R_V$ and $R_\varepsilon$ defined in (3.8) are now shown to vanish asymptotically. At first, consider $R_V$ defined in (3.8), which consists of the terms involving the unit root components $V_t$. Several results are given, of which the strongest result for $R_V$ can only be established for certain values of the parameters.

**Lemma 3.5.** *Suppose $A_k = 0$ and that the Assumptions $(2.1)$, $(2.3)$, $(2.5)$ are satisfied with $\lambda > 2$. Then*

- (i) $R_V \stackrel{a.s.}{=} \mathrm{O}(\log T)$,
- (ii) $R_V = \mathrm{o_P}(1)$ *if also Assumption $(2.2)$ holds,*
- (iii) $R_V \stackrel{a.s.}{=} \mathrm{O}\{(\log\log T)^{1/2}(\log T)^{1/2}\}$ *if $(A)$ $\mathbf{D}$ and $\mathbf{V}$ have no common eigenvalues,*
- (iv) $R_V \stackrel{a.s.}{=} \mathrm{o}(1)$ *if $(B)$ $\dim \mathbf{D} = 0$ and $\mathbf{V} = 1$,*
- (v) $R_V = 0$ *if $(C)$ $\dim \mathbf{V} = 0$.*

*Proof of Lemma 3.5.* (i) This follows since $Q(V_{t-j}|D_t) \stackrel{a.s.}{=} \mathrm{O}(\log T)$ according to Nielsen ([17], Theorem 2.4).

(ii) The type of argument for (3.10) in the proof of Lemma 3.4 can be used. Replacing $W$ with $V$ throughout, the asymptotic properties of $a_T, b_T, c_T, d_T$ have to be explored. For $b_T, c_T$ the argument is the same so, for all $\eta > 0$,

$$b_T \stackrel{a.s.}{=} \mathrm{O}(T), \qquad c_T \stackrel{a.s.}{=} \mathrm{o}(T^{1/2+\eta}),$$

whereas using Nielsen ([17], Theorems 2.4) for $a_T$ and the techniques of Chan and Wei [5] for $d_T$ shows, for all $\eta > 0$,

$$a_T \stackrel{a.s.}{=} \mathrm{o}\left(T^\eta\right), \qquad d_T^{-1} = \mathrm{o_P}(T^{-1-4\eta}),$$

so $f_T = \mathrm{o_P}(T^{-4\eta})$. Since $Q(V_{t-j}|D_t) \stackrel{a.s.}{=} \mathrm{O}(\log T)$ as established in (*i*) the desired result follows by expanding $Q(V_{t-1}|D_t)$ in terms of $Q(V_{t-2}|D_t)$.

(iii) Define the vector $S_{t-1} = (V_{t-1}', D_t')'$. By partitioned inversion it holds

$$Q\left(S_{t-1}\right) = Q\left(V_{t-1}|D_t\right) + Q\left(D_t\right).$$

By an invariance argument $D_t$ can be replaced by $D_{t-j}$ and thus it follows

$$R_V = Q\left(V_{t-2}|D_t\right) - Q\left(V_{t-1}|D_t\right) = Q\left(S_{t-2}\right) - Q\left(S_{t-1}\right).$$

Due to (2.4) and (2.7) the process $S_{t-1}$ satisfies $S_t = \mathbf{S}S_{t-1} + e_{S,t}$ for a matrix $\mathbf{S}$ with eigenvalues of length one and $e_{S\,t} = (e'_{V,t}, 0')'$. It then follows that

$$\sum_{t=1}^{T} \varepsilon_t S'_{t-1} = \sum_{t=1}^{T} \varepsilon_t \left( S'_{t-2}\mathbf{S}' + e'_{S,t-1} \right).$$

Inserting this expression into $Q(S_{t-1})$ shows

$$Q\left(S_{t-1}\right) = \sum_{t=1}^{T} \varepsilon_t S'_{t-1} \left( \sum_{t=1}^{T} S^{\otimes 2}_{t-1} \right)^{-1} \sum_{t=1}^{T} S_{t-1} \varepsilon'_t = Q_A + Q_B + Q_C + Q'_C,$$

where

$$Q_A = Q_1 Q_2 Q'_1, \quad Q_B = Q_4 Q_3 Q'_3 Q'_4, \quad Q_C = Q_1 Q_2^{1/2} Q'_3 Q'_4,$$

are defined in terms of the statistics

$$Q_1 = \sum_{t=1}^{T} \varepsilon_t e'_{S,t-1}, \quad Q_2 = \left( \sum_{t=1}^{T} S^{\otimes 2}_{t-1} \right)^{-1},$$

$$Q_3 = \left( \sum_{t=1}^{T} S^{\otimes 2}_{t-2} \right)^{1/2} \mathbf{S} \left( \sum_{t=1}^{T} S^{\otimes 2}_{t-1} \right)^{-1/2}, \quad Q_4 = \sum_{t=1}^{T} \varepsilon_t S'_{t-2} \left( \sum_{t=1}^{T} S^{\otimes 2}_{t-2} \right)^{-1/2}.$$

The orders of magnitude of these follow from a series of results in Nielsen [17]. Theorem 6.1 and Lemma 6.3 imply $Q_1 \stackrel{a.s.}{=} \mathrm{O}\{(T \log\log T)^{1/2}\}$. Theorem 8.3 shows $Q_2 \stackrel{a.s.}{=} \mathrm{O}(T^{-1})$ when $\mathbf{D}$ and $\mathbf{V}$ have no common eigenvalues. Lemma 8.7$(ii)$ shows $Q_3^{\otimes 2} - I \stackrel{a.s.}{=} \mathrm{O}\{T^{-1/2}(\log T)^{1/2}\}$. Theorem 2.4 shows $Q_4 \stackrel{a.s.}{=} \mathrm{O}\{(\log T)^{1/2}\}$. Noting that $Q(S_{t-2}) = Q_4 Q'_4$ this in turn implies

$$Q_A = \mathrm{O}(\log\log T), \qquad Q_B = Q(S_{t-2}) + \mathrm{O}\{T^{-1/2}(\log T)^{3/2}\},$$
$$Q_C = \mathrm{O}\{(\log\log T)^{1/2}(\log T)^{1/2}\},$$

and the desired result follows.

(iv) Donsker and Varadhan's [6] Law of the Iterated Logarithm for the integrated squared Brownian motion states

$$\liminf_{T \to \infty} \frac{\log\log T}{T^2} \int_0^T B_u^2 du \stackrel{a.s.}{=} \frac{1}{4}.$$

Now use *either* the argument in (ii) with $d_T^{-1} \stackrel{a.s.}{=} \mathrm{O}(T^{-2} \log\log T)$ *or* the argument in (iii) with $Q_2 \stackrel{a.s.}{=} \mathrm{O}(T^{-2} \log\log T)$ so $Q_A, Q_B, Q_C$ are all o(1).

(v) This follows by construction. $\qquad\square$

Now, consider $R_\varepsilon$ defined in (3.8). By showing that this vanishes it follows that the regressors can be excluded asymptotically in the term involving the lagged innovations $\varepsilon_{t-1}$. A fourth order moment condition is now needed in Assumption (2.1).

**Lemma 3.6.** *Suppose $A_k = 0$ and that the Assumptions $(2.1)$, $(2.3)$, $(2.5)$ are satisfied, now with $\lambda > 4$. Then*

$$R_\varepsilon = Q\left(\varepsilon_{t-1}|\mathbf{X}_{t-2}, D_t\right) - Q\left(\varepsilon_{t-1}\right) \stackrel{a.s.}{=} \mathrm{o}\left(1\right).$$

*Proof of Lemma 3.6.* Define the vector $S_t = (\mathbf{X}'_{t-2}, D'_t)'$. According to Nielsen ([17], Theorem 2.4) it holds that, for any $\eta > 0$, the terms

$$(3.11) \qquad \left(\sum_{t=1}^T S_t S'_t\right)^{-1/2} \sum_{t=1}^T S_t \varepsilon'_t, \qquad \left(\sum_{t=1}^T S_t S'_t\right)^{-1/2} \sum_{t=1}^T S_t \varepsilon'_{t-1}$$

are $\mathrm{o}(T^{1/4-\eta})$ when indeed $\lambda > 4$. It then holds that

$$\sum_{t=1}^T \varepsilon_t \varepsilon'_{t-1} - \sum_{t=1}^T \varepsilon_t S'_t \left(\sum_{t=1}^T S_t S'_t\right)^{-1} \sum_{t=1}^T S_t \varepsilon'_{t-1} \overset{a.s.}{=} \sum_{t=1}^T \varepsilon_t \varepsilon'_{t-1} + \mathrm{o}(T^{1/2-\eta}),$$

$$\sum_{t=1}^T \varepsilon_{t-1} \varepsilon'_{t-1} - \sum_{t=1}^T \varepsilon_{t-1} S'_t \left(\sum_{t=1}^T S_t S'_t\right)^{-1} \sum_{t=1}^T S_t \varepsilon'_{t-1} \overset{a.s.}{=} \sum_{t=1}^T \varepsilon_{t-1} \varepsilon'_{t-1} + \mathrm{o}(T^{1-\eta}),$$

where the requirement $\lambda > 4$ is only needed in the first case. Theorems 2.5, 6.1 of the above paper show $T^{-1} \sum_{t=1}^T \varepsilon_{t-1} \varepsilon'_{t-1}$ has positive definite limit points while $\sum_{t=1}^T \varepsilon_t \varepsilon'_{t-1} (\sum_{t=1}^T \varepsilon_{t-1} \varepsilon'_{t-1})^{-1/2} = \mathrm{o}(T^\eta)$. Combine these results. □

**Remark 3.7.** In Lemma 3.6 a fourth moment condition comes in through the requirement that $\lambda > 4$ in Assumption (2.1). This can be relaxed to $\lambda > 2$ under one of two alternative assumptions.

(I,a) If $\dim \mathbf{W} = 0$ then the terms in (3.11) are $\mathrm{o}(T^\eta)$, see Nielsen ([17], Theorem 2.4), and the main result holds.

(I,b) If $\dim \mathbf{W} > 0$ but the innovations $\varepsilon_t$ are independently, identically distributed then terms of the type $(\sum_{t=1}^T W_{t-1} W'_{t-1})^{-1/2} \sum_{t=1}^T W_{t-1} \varepsilon'_t$ converge in distribution, see Anderson [1] and the result of the Theorem holds, albeit only in probability.

### 3.3.4. *The leading term of $\hat{\Omega}_{k-1} - \hat{\Omega}_k$*

First the order of magnitude the leading term in (3.7) is established in an almost sure sense. This can be done under weak moment conditions. Subsequently the distribution of the leading term is investigated.

**Lemma 3.8.** *Suppose $A_k = 0$ and that the Assumptions (2.1), (2.3) are satisfied with $\lambda > 2$. Define $E_T = T^{-1} \sum_{t=1}^T \varepsilon_t \varepsilon'_t$. Then*
$\limsup_{T \to \infty} (2 \log \log T)^{-1} \mathrm{tr}[\{Q(\varepsilon_{t-1}) + Q(\tilde{U}_{t-2}) - Q(\tilde{U}_{t-1})\} E_T^{-1}] \overset{a.s.}{=} \mathrm{O}(1).$

*Proof of Lemma 3.8.* This follows by noting that the sequence $\hat{\Omega}_{k-1}^{-1}$ is relatively compact with positive definite limiting points due to Lemma 3.1 and Lai and Wei ([15], Theorem 2) and otherwise following the argument in the proof of Pötscher ([21], Theorem 3.4). □

When it comes to analysing the distribution of the leading term in (3.7) it is convenient to show that it can be written as a single quadratic form $Q(Y_{t-1})$ for some process $Y_{t-1}$. This argument requires two steps, of which the first is concerned with the convergence properties of $T^{-1} \sum_{t=1}^T \tilde{U}_{t-1} \tilde{U}'_{t-1}$. As the argument involves a variance matrix, the Assumption (2.2) is now called upon.

**Lemma 3.9.** *Suppose $A_k = 0$ and that the Assumptions $(2.1), (2.2)$ are satisfied with $\lambda > 2$. Let $M_U$ be the matrix defined by $e_{U,t} = M_U \varepsilon_t$ in $(2.7)$ and define*

$$F = \sum_{t=0}^{\infty} \mathbf{U}^t M_U \Omega M_U' (\mathbf{U}^t)'.$$

*Then for all $\zeta < \min(1 - 2/\lambda, 1/2)$ it holds*

$$\frac{1}{T} \sum_{t=1}^{T} \tilde{U}_t \tilde{U}_t' \overset{a.s.}{=} F + \mathrm{o}\left(T^{-\zeta}\right).$$

*Proof of Lemma 3.9.* Following the proof of Lai and Wei ([15], Theorem 2), the equation $(2.8)$ shows

$$\sum_{t=1}^{T} \tilde{U}_t \tilde{U}_t' \overset{a.s.}{=} \mathbf{U} \left( \sum_{t=1}^{T} \tilde{U}_t \tilde{U}_t' - \tilde{U}_T \tilde{U}_T' + \tilde{U}_0 \tilde{U}_0' \right) \mathbf{U}'$$

$$+ M_U \sum_{t=1}^{T} \varepsilon_t \varepsilon_t' M_U' + \mathrm{O}\left( \sum_{t=1}^{T} \tilde{U}_{t-1} \varepsilon_t' \right).$$

Due to Nielsen ([17], Theorems 2.4, 5.1, Example 6.5) both $\sum_{t=1}^{T} \tilde{U}_{t-1} \varepsilon_t'$ and $\tilde{U}_T \tilde{U}_T'$ are $\mathrm{o}(T^{1-\zeta})$. Note that Assumption $(2.5)$ is not needed as $\tilde{U}_t$ does not involve deterministic terms. Denoting $F_T = T^{-1} \sum_{t=1}^{T} \tilde{U}_t \tilde{U}_t'$ it follows from Lemma 3.1 that

$$F_T - \mathbf{U} F_T \mathbf{U}' \overset{a.s.}{=} M_U \Omega M_U' + \mathrm{o}(T^{-\zeta}).$$

This equation has a unique solution $F_T = \sum_{t=0}^{\infty} \mathbf{U}^t \{ M_U E M_U' + \mathrm{o}(T^{-\zeta}) \} (\mathbf{U}^t)'$, see Anderson and Moore ([2], p. 336), which in turn equals $F + \mathrm{o}(T^{-\zeta})$ since the maximal eigenvalue of $\mathbf{U}\mathbf{U}'$ is less than one. $\qquad\square$

The leading term in $(3.7)$ is now written as a single quadratic form $Q(Y_{t-1})$.

**Lemma 3.10.** *Suppose $A_k = 0$ and that the Assumptions $(2.1), (2.2)$ are satisfied with $\lambda > 2$. Then there exists an $\{(p + \dim U) \times p\}$-matrix $C$ with full column rank so*

$$Q\left(\varepsilon_{t-1}\right) + Q(\tilde{U}_{t-2}) - Q(\tilde{U}_{t-1}) \overset{a.s.}{=} Q\left(Y_{t-1}\right) + \mathrm{o}\left(1\right),$$

*where $Y_t$ is the process $C'(\varepsilon_t', U_{t-1}')'$.*

*Proof of Lemma 3.10.* The idea is to exploit that the asymptotic covariance for $Z_{t-1} = (\tilde{U}_{t-2}', \varepsilon_{t-1}')'$ is diagonal with elements $F, \Omega$. By the above Lemmas 3.1, 3.9 then, for some $\eta > 0$,

$$Q\left(\varepsilon_{t-1}\right) + Q(\tilde{U}_{t-2})$$

$$(3.12) \quad = \sum_{t=1}^{T} \varepsilon_t \begin{pmatrix} \tilde{U}_{t-2} \\ \varepsilon_{t-1} \end{pmatrix}' \left\{ \sum_{t=1}^{T} \begin{pmatrix} \tilde{U}_{t-2} \tilde{U}_{t-2}' & 0 \\ 0 & \varepsilon_{t-1} \varepsilon_{t-1}' \end{pmatrix} \right\}^{-1} \sum_{t=1}^{T} \begin{pmatrix} \tilde{U}_{t-2} \\ \varepsilon_{t-1} \end{pmatrix} \varepsilon_t'$$

$$\overset{a.s.}{=} \frac{1}{T} \sum_{t=1}^{T} \varepsilon_t \begin{pmatrix} \tilde{U}_{t-2} \\ \varepsilon_{t-1} \end{pmatrix}' \begin{pmatrix} F & 0 \\ 0 & \Omega \end{pmatrix}^{-1} \sum_{t=1}^{T} \begin{pmatrix} \tilde{U}_{t-2} \\ \varepsilon_{t-1} \end{pmatrix} \varepsilon_t' \{ 1 + \mathrm{o}(T^{-\eta}) \}$$

As discussed in Section 2 then $\tilde{U}_{t-1} = \mathbf{U} \tilde{U}_{t-2} + M_U \varepsilon_{t-1}$ for some matrix $M_U$ with full column rank. In particular $\tilde{U}_{t-1} = C_\perp' (\tilde{U}_{t-2}', \varepsilon_{t-1}')'$ where the $\{(p + \dim U) \times$ dim $U$}-matrix $C_\perp = (\mathbf{U}, M_U)'$ has full column rank. Therefore a $\{(p+\dim U) \times p\}$-matrix $C$ can be chosen with full column rank so the matrix $(C, C_\perp)$ is regular and

$$C' \begin{pmatrix} F & 0 \\ 0 & \Omega \end{pmatrix} C_\perp = 0.$$

The sequences $T^{-1} \sum_{t=1}^{T} \tilde{U}_{t-1} \tilde{U}'_{t-1}$ and $T^{-1} \sum_{t=1}^{T} \tilde{U}_{t-2} \tilde{U}'_{t-2}$ will have the same limit, $F$, while $T^{-1} \sum_{t=1}^{T} Y_{t-1} Y'_{t-1}$ will converge to a positive definite matrix $G$. It then holds

$$\begin{pmatrix} C'_\perp \\ C' \end{pmatrix} \begin{pmatrix} F & 0 \\ 0 & \Omega \end{pmatrix} (C_\perp, C) = \begin{pmatrix} F & 0 \\ 0 & G \end{pmatrix}.$$

Pre- and post-multiplying the middle matrix in (3.12) with $(C_\perp, C)(C_\perp, C)^{-1}$ and its transpose then implies

$$Q(\varepsilon_{t-1}) + Q(\tilde{U}_{t-2}) \stackrel{a.s.}{=} \left\{ Q(\tilde{U}_{t-1}) + Q(Y_{t-1}) \right\} \{1 + \mathrm{o}(T^{-\eta})\}.$$

Theorem 2.4 of Nielsen (2005) implies $Q(\tilde{U}_{t-1})$ and $Q(Y_{t-1})$ are $\mathrm{o}(T^\eta)$, which gives the desired result. □

The asymptotic distribution of the leading term $Q(Y_{t-1})$ now follows.

**Lemma 3.11.** *Suppose $A_k = 0$ and that the Assumptions (2.1), (2.2), (2.3) are satisfied with $\lambda > 4$. Then*

(i) $1 \le \limsup_{T \to \infty} (2 \log \log T)^{-1} \mathrm{tr}\{Q(Y_{t-1})\Omega^{-1}\} \le p^2$ *a.s.*

(ii) $\mathrm{tr}\{Q(Y_{t-1})\Omega^{-1}\} \xrightarrow{\mathrm{D}} \chi^2(p^2)$.

*Proof of Lemma 3.11.* (i) This follows from the Law of Iterated Logarithms by Heyde and Scott ([12], Corollary 2) and Hannan ([9], p. 1076-1077). See Quinn [22] for details.

(ii) This follows from Brown and Eagleson's [4] Central Limit Theorem. This requires existence of second moments of $\varepsilon_t Y_{t-1}$. □

**Remark 3.12.** The proof of Lemma 3.11 actually only requires the existence of fourth moments, which is slightly weaker than the stated condition of $\lambda > 4$ in Assumption (2.1). In Lemma 3.11(ii) this can be relaxed to a second moment condition if for instance:
(II) the innovations $\varepsilon_t$ are independent.

### 3.4. Proofs of results for likelihood ratio test statistics

*Proof of Theorem 2.1.* Consider the formula (3.1). The term $\hat{\Omega}_{k-1}$ was dealt with in Lemma 3.1. As for the term $T(\hat{\Omega}_{k-1} - \hat{\Omega}_k)$ consider two cases.

When $k = 1$ then $T(\hat{\Omega}_{k-1} - \hat{\Omega}_k) = Q(\varepsilon_{t-1})$.

When $k > 1$ apply the expansion in Lemma 3.4. The term $R_V$ vanishes due to Lemma 3.5(ii) when Assumption (2.2) is satisfied. The term $R_\varepsilon$ vanishes due to Lemma 3.6 when $\lambda > 4$ in Assumption (2.1). Due to Lemma 3.10 the leading term is now $Q(Y_{t-1})$, provided Assumption (2.2) holds.

For any $k$ the desired $\chi^2$-distribution now arises from Lemma 3.11(ii) provided Assumptions (2.2), (2.1) are satisfied with $\lambda > 4$. □

*Proof of Theorem 2.3.* Note first that $T(\hat{\Omega}_{k-1} - \hat{\Omega}_{k+m-1})$ can be written as $Q(\tilde{X}_{t-k}^{t-k-m+1}|\mathbf{X}_{t-1}, D_t)$ where $\tilde{X}_{t-a}^{t-b} = (X'_{t-a}, \ldots, X'_{t-b})'$. Consider now the proof of the decomposition in Lemma 3.2. Using first (3.4) and then (3.5) repeatedly it is seen that

$$T\left(\hat{\Omega}_{k-1} - \hat{\Omega}_{k+m-1}\right) = Q\left(\left.\begin{matrix}\mathbf{X}_{t-1}\\ \tilde{X}_{t-k}^{t-k-m+1}\end{matrix}\right| D_t\right) - Q\left(\mathbf{X}_{t-1}| D_t\right)$$

$$= \sum_{j=1}^{m} Q\left(\varepsilon_{t-j}| X_{t-j-1}^{t-m}, \mathbf{X}_{t-m-1}, D_t\right)$$

$$+ Q\left(\mathbf{X}_{t-m-1}| D_t\right) - Q\left(\mathbf{X}_{t-1}| D_t\right).$$

As in the proof of Theorem 2.1 the Lemmas 3.4, 3.5($ii$), 3.6 show that the leading terms reduce to

$$T\left(\hat{\Omega}_{k-1} - \hat{\Omega}_{k+m-1}\right) = \sum_{j=1}^{m} Q\left(\varepsilon_{t-j}\right) + Q\left(\tilde{U}_{t-m-1}\right) - Q\left(\tilde{U}_{t-1}\right) + o_{\mathsf{P}}\left(1\right),$$

when $k_0 < k$. A slight generalisation of Lemma 3.10 is needed, using that the asymptotic covariance for $Z_{t-1} = (\tilde{U}'_{t-m-1}, \varepsilon'_{t-1}, \ldots, \varepsilon'_{t-m})'$ is diagonal with elements $F, \Omega, \ldots, \Omega$. A $\{(mp + \dim U) \times mp\}$-matrix $C$ can then be found giving rise to a process $Y_{t-1} = C'Z_{t-1}$. The argument is completed using a Central Limit Theorem as in the proof of Lemma 3.11($ii$).  □

## 3.5. Proofs of results for information criteria

*Proof of Theorem 2.4.* Consider $j < k_0$. The condition $f(T) = o(T)$ implies

$$\Phi_j - \Phi_{k_0} = \log \det\{I + (\hat{\Omega}_j - \hat{\Omega}_{k_0})\hat{\Omega}_{k_0}^{-1}\} + o\left(1\right).$$

Lemma 3.1 shows that $\hat{\Omega}_{k_0} \overset{a.s.}{\to} \Omega$, so it suffices that $\liminf_{T\to\infty} \lambda_{\max}(\hat{\Omega}_j - \hat{\Omega}_{k_0})$ is positive. Defining $\mathbf{Y}_t = (X'_{t-1}, \ldots, X'_{t-j+1})'$ and $\mathbf{Z}_t = (X'_{t-j}, \ldots, X'_{t-k_0})'$ it holds

$$\hat{\Omega}_j - \hat{\Omega}_{k_0} = \left[T^{-1/2}\sum_{t=1}^{T} X_t \left(\mathbf{Z}_{t-1}|\mathbf{Y}_{t-1}, D_t\right)'\left\{\sum_{t=1}^{T}\left(\mathbf{Z}_{t-1}|\mathbf{Y}_{t-1}, D_t\right)^{\otimes 2}\right\}^{-1/2}\right]^{\otimes 2}.$$

Define $\mathbf{A_y} = A_1, \ldots, A_j$ and $\mathbf{A_z} = A_{j+1}, \ldots, A_{k_0}$ noting that $A_{k_0} \neq 0$. Then it holds $X_t = \mathbf{A_y}\mathbf{Y}_t + \mathbf{A_z}\mathbf{Z}_t + \mu D_t + \varepsilon_t$. Therefore $\hat{\Omega}_j - \hat{\Omega}_{k_0}$ equals

$$T^{-1/2}\sum_{t=1}^{T} \varepsilon_t \left(\mathbf{Z}_{t-1}|\mathbf{Y}_{t-1}, D_t\right)'\left\{\sum_{t=1}^{T}\left(\mathbf{Z}_{t-1}|\mathbf{Y}_{t-1}, D_t\right)^{\otimes 2}\right\}^{-1/2}$$

$$+\mathbf{A_z}\left\{T^{-1}\sum_{t=1}^{T}\left(\mathbf{Z}_{t-1}|\mathbf{Y}_{t-1}, D_t\right)^{\otimes 2}\right\}^{1/2}.$$

The first term is of order o(1) *a.s.* by Nielsen ([17], Theorem 2.4). As for the second term it holds that $\liminf_{T\to\infty} \lambda_{\min}\{T^{-1}\sum_{t=1}^{T}(\mathbf{X}_{t-1}|D_t)^{\otimes 2}\} > 0$ *a.s.* according to Nielsen ([17], Corollary 9.5). As a consequence the limit points of $T^{-1}\sum_{t=1}^{T}(\mathbf{Z}_{t-1}|\mathbf{Y}_{t-1}, D_t)^{\otimes 2}$ are positive definite. Since $\mathbf{A_z} \neq 0$ then $\liminf_{T\to\infty} \lambda_{\min}(\hat{\Omega}_j - \hat{\Omega}_{k_0}) > 0$ and therefore $\liminf_{T\to\infty} \hat{k} \geq k_0$ *a.s.*  □

*Proof of Theorem 2.5.* Consider now $k_0 < j \le K$. It then holds

$$\Phi_{j+1} - \Phi_j = \log\det(\hat\Omega_{j+1}\hat\Omega_j^{-1}) + T^{-1}f\,(T)$$
$$= \log\det\{I_p - (\hat\Omega_j - \hat\Omega_{j+1})\hat\Omega_j^{-1}\} + T^{-1}f\,(T)\,.$$

A Taylor expansion shows

$$\Phi_{j+1} - \Phi_j \stackrel{a.s.}{=} -\mathrm{tr}\{(\hat\Omega_j - \hat\Omega_{j+1})\hat\Omega_j^{-1}\} + T^{-1}f\,(T) + \mathrm{o}[\{(\hat\Omega_j - \hat\Omega_{j+1})\hat\Omega_j^{-1}\}^2].$$

Lemma 3.1 shows that $\hat\Omega_j$ is consistent, while Lemma 3.4 gives the expansion

$$T(\hat\Omega_{j-1} - \hat\Omega_j) \stackrel{a.s.}{=} Q\,(\varepsilon_{t-1}) + Q(\tilde U_{t-2}) - Q(\tilde U_{t-1}) + R_\varepsilon + R_V + \mathrm{o}\,(1)\,.$$

To complete the proof it has to be shown that $\Phi_{j+1} - \Phi_j$ has a positive limiting value. This holds if $T(\hat\Omega_{j-1} - \hat\Omega_j) = \mathrm{o}\{g(T)\}$ for some function $g(T)$ so $f(T)/g(T) \to \infty$.

(i) The term $R_V$ vanishes due to Lemma 3.5(ii) when Assumption (2.2) is satisfied. The term $R_\varepsilon$ vanishes due to Lemma 3.6 when $\lambda > 4$ in Assumption (2.1). Due to Lemma 3.10 the leading term is $Q(Y_{t-1})$, provided Assumption (2.2) holds. This is $\mathrm{O}_\mathsf{P}(1)$ by Lemma 3.11($ii$) provided Assumptions (2.1), (2.2) are satisfied with $\lambda > 4$.

(ii) The term $R_V$ is $\mathrm{O}(\log T)$ due to Lemma 3.5(i). The term $R_\varepsilon$ vanishes due to Lemma 3.6 when $\lambda > 4$ in Assumption (2.1). Due to Lemma 3.8 the leading term is $\mathrm{O}(\log\log T)$.

(iii) Under $(A)$ that $\mathbf{V}$ and $\mathbf{D}$ have no common eigenvalues then $R_V$ is $\mathrm{O}\{(\log T)^{1/2}(\log\log T)^{1/2}\}$ due to Lemma 3.5($iii$). The argument of ($ii$) can then be followed.

(iv) Under $(B)$ that $\dim\mathbf{D} = 0$ with $\mathbf{V} = 1$ then $R_V$ is $\mathrm{o}(1)$ due to Lemma 3.5($iv$), whereas under $(C)$ that $\dim\mathbf{V} = 0$ then $R_V = 0$. while it is $\mathrm{o}(1)$ under $(B)$ $\dim\mathbf{D} = 0$ with $\mathbf{V} = 1$ due to Lemma 3.5($iv$). The argument of ($ii$) can then be followed.

(v) The terms $R_V$ and $R_\varepsilon$ vanish as in (iv). As in $(i)$ the leading term is $Q(Y_{t-1})$ by Lemma 3.10 provided Assumption (2.2) holds. This is of the desired order of magnitude by Lemma 3.11($i$) provided Assumptions (2.2), (2.1) are satisfied with $\lambda > 4$. □

**Remark 3.13.** The condition $\lambda > 4$ in Theorem 2.5 can be relaxed as follows.

(i) It is used first in Lemma 3.6 and can be relaxed under (I,a) or (I,b) as this is a result holding in probability, see Remark 3.7. It is used secondly in Lemma 3.11(ii) and can be relaxed under (II), see Remark 3.12.

(ii), (iii), (iv) It is only used in Lemma 3.6 and can only be relaxed under (I,a) as this is a result holding almost surely, see Remark 3.7.

(v) It is indeed required in Lemma 3.11(i).

### 3.6. *Proof of results for residual based tests*

*Proof.* It suffices to show how the residual based test statistics relate to the likelihood ratio test statistics.

In the joint test the squared sample multiple correlation $R^2$ of $\hat\varepsilon_t$ and the vector $Z_{t-1} = (\hat\varepsilon'_{t-1}, \ldots, \hat\varepsilon'_{t-m}, \mathbf{X}'_{t-1}, D'_t)'$ is considered, recalling that $\mathbf{X}_{t-1}$ is defined as $(X'_{t-1}, \ldots, X'_{t-k+1})'$. The key to the result is that

$$\hat\varepsilon_{t-j} = X_{t-j} - \hat{\mathbf{B}}\mathbf{X}_{t-j-1} - \hat\mu\mathbf{D}^{-j-1}D_t,$$

where $\hat{\mathbf{B}}, \hat{\mu}$ are least squares estimators based on (1.1) for the full sample $t = 1, \ldots, T$. Due to the inclusion of $\mathbf{X}_{t-1}$ as regressor it follows that $Z_{t-1} = N\tilde{Z}_{t-1}$ where $\tilde{Z}_{t-1} = (X_{t-1}, \ldots, X_{t-k-m+1}, D_t)$ and the square matrix $N$ is based on $\hat{\mathbf{B}}, \hat{\mu}$ and is invertible with probability one. By the invariance of sample multiple correlations to linear transformations then $R^2$ can be computed from $\hat{\varepsilon}_t$ and $\tilde{Z}_{t-1}$. By the same type of manipulation as in Lemma 3.2 it follows that

$$\hat{Q}\left(\tilde{Z}_{t-1}\right) = \sum_{t=m+1}^{T} \hat{\varepsilon}_t \tilde{Z}'_{t-1} \left( \sum_{t=m+1}^{T} \tilde{Z}_{t-1}^{\otimes 2} \right)^{-1} \sum_{t=m+1}^{T} \tilde{Z}_{t-1} \hat{\varepsilon}'_t$$

can be written as

$$(3.13) \qquad \hat{Q}\left(\tilde{Z}_{t-1}\right) = \hat{Q}\left(X_{t-k}, \ldots, X_{t-k-m+1} | \mathbf{X}_{t-1}, D_t\right) + \hat{Q}\left(\mathbf{X}_{t-1}, D_t\right).$$

Since the first term in (3.13) includes the regressors $\mathbf{X}_{t-1}, D_t$ then $\hat{\varepsilon}_t$ can be replaced by $\varepsilon_t$. Thus, apart from starting the regression at $t = m + 1$ instead of $t = 1$ this term is the same as $Q(X_{t-k}, \ldots, X_{t-k-m+1} | \mathbf{X}_{t-1}, D_t)$. It therefore has the same asymptotic properties as $T(\hat{\Omega}_{k-1} - \hat{\Omega}_{k+m-1})$, which was studied in the proof of Theorem 2.3.

The second term in (3.13) vanishes asymptotically. This is because the residuals $\hat{\varepsilon}_t$ are orthogonal to $\mathbf{X}_{t-1}, D_t$ when evaluated over $t = 1, \ldots, T$. A tedious analysis shows that this orthogonality holds asymptotically when evaluated over $t = m + 1, \ldots, T$.

For the marginal test the argument is the same. The main difference is that the residuals are now

$$\hat{\varepsilon}_{t-j,\mathrm{marg}} = X_{t-j,1} - \hat{\mathbf{B}}_{\mathrm{marg}}\mathbf{X}_{t-j-1} - \hat{\mu}_{\mathrm{marg}}\mathbf{D}^{-j-1}D_t.$$

Once again the inclusion of $\mathbf{X}_{t-1}$ as regressor implies that the vector $Z_{t-1,\mathrm{marg}}$ defined as $(\hat{\varepsilon}'_{t-1}, \ldots, \hat{\varepsilon}'_{t-m}, \mathbf{X}'_{t-1}, D'_t)'$ can be replaced by the above $\tilde{Z}_{t-1}$. So the statistic $\hat{Q}(\tilde{Z}_{t-1})$ is replaced by a statistic based on $\hat{\varepsilon}_{t,\mathrm{marg}}$, but the same $\tilde{Z}_{t-1}$.

For the conditional test the residuals are of the type

$$\hat{\varepsilon}_{t-j,\mathrm{cond}} = X_{t-j,1} - \hat{\mathbf{B}}_{\mathrm{cond}}\mathbf{X}_{t-j-1} - \hat{\mu}_{\mathrm{cond}}\mathbf{D}^{-j-1}D_t - \hat{\omega}X_{t-j,2}.$$

The same argument applies as for the marginal test. $\qquad\square$

## Acknowledgments

Comments from the referee are gratefully acknowledged.